\documentclass[11pt]{amsart}

\usepackage{amsmath, amssymb, amsthm, amscd}
\usepackage[top=1 in, bottom= 1 in, left= 1 in, right= 1 in]{geometry}
\usepackage{cite}
\usepackage{graphicx}
\usepackage{comment}

\newcommand{\ds}{\displaystyle}

\theoremstyle{plain}
\newtheorem{theorem}{Theorem}[section]
\newtheorem{lemma}[theorem]{Lemma}
\newtheorem{proposition}[theorem]{Proposition}

\theoremstyle{definition}
\newtheorem*{definition}{Definition}
\newtheorem{example}[theorem]{Example}
\newtheorem{remark}[theorem]{Remark}
\newtheorem*{note}{\underline{Note}}
\newtheorem*{acknowledgments}{Acknowledgments}

\numberwithin{equation}{subsection}
\numberwithin{theorem}{section}

\newcommand{\frattini}[1]{[#1, #1]#1^{p}}

\begin {document}

\title[$RG$ \& $RG_{p}$ of Amalgams \& HNN Extensions]{Rank Gradient \& $p$-Gradient of Amalgamated Free Products \& HNN Extensions}
\author{Nathaniel Pappas\\University of Virginia}
\maketitle

\begin{abstract}
We calculate the rank gradient and $p$-gradient of free products with amalgamation over an amenable subgroup and HNN extensions with an amenable associated subgroup.  The notion of cost is used to compute the rank gradient of amalgamated free products and HNN extensions.  For the $p$-gradient the Kurosh subgroup theorems for amalgamated free products and HNN extensions will be used.    
\end{abstract}

\section{Introduction}
The rank gradient and $p$-gradient are two group invariants which originated in topology.  Mark Lackenby first introduced the rank gradient \cite{LacExpanders} and $p$-gradient \cite{LackenbyTau} as means to study $3$-manifold groups.  However, these group invariants have been gaining interest among group theorists. Both invariants are difficult to compute and for the majority of classes of groups where the rank gradient has been calculated its value is zero. We add to the few results on computing rank gradient and $p$-gradient by giving formulas for the rank gradient and $p$-gradient of free products with amalgamation over an amenable subgroup and HNN extensions with an amenable associated subgroup.  

Let $\Gamma$ be a finitely generated group and let $d(\Gamma)$ denote the minimal number of generators of $\Gamma$.  A set of subgroups $\{H_{n}\}$ of $\Gamma$ is called a \textit{lattice} if it is closed under finite intersections.  In particular any descending chain of subgroups is a lattice. The \textit{rank gradient} relative to a lattice $\{H_{n}\}$ of finite index subgroups is defined as
\[RG(\Gamma, \{H_{n}\}) = \inf_{n} \frac{d(H_{n}) - 1}{[\Gamma : H_{n}]}.\]

For a prime number $p$ the related notion of $p$-gradient, $RG_{p}(\Gamma, \{H_{n}\})$, is defined similarly by replacing $d(H)$ with $d_{p}(H)= d(H/[H,H]H^{p})$ and requiring the subgroups to be normal of $p$-power index.  One can also define the \textit{absolute rank gradient}, $RG(\Gamma)$, (resp. $p$-gradient, $RG_{p}(\Gamma)$) where the infimum is taken over all finite index subgroups (resp. normal subgroups of $p$-power index).
  
\begin{remark}
All of the results given below are stated for the rank gradient but the analogous results hold for $p$-gradient for any prime $p$, and are explicitly stated in Section~\ref{p grad sec}.  
\end{remark}

Abert, Jaikin-Zapirain, and Nikolov in \cite{AJN} computed the rank gradient of a free product of residually finite groups relative to a descending chain of normal subgroups using Bass-Serre theory.  

\begin{theorem}[Abert, Jaikin-Zapirain, and Nikolov]
\label{RG of free prod rel chain}
Let $\Gamma_{1}$ and $\Gamma_{2}$ be finitely generated and residually finite.  Let $\{H_{n}\}$ be a normal chain of finite index subgroups in $\Gamma = \Gamma_{1} \ast \Gamma_{2}$.  Then
\[RG(\Gamma, \{H_{n}\}) = RG(\Gamma_{1}, \{\Gamma_{1} \cap H_{n}\}) + RG(\Gamma_{2}, \{\Gamma_{2} \cap H_{n}\}) +1.\]
\end{theorem}

The difficulty with obtaining similar rank results for free products with amalgamation or HNN extensions is getting a lower bound on the minimal number of generators of a finite index subgroup due to the lack of a Grushko-Neumann like result in these cases. To get around this issue, we use the theory of \textit{cost}. Rank gradient is closely related to cost as well as L$^{2}$-Betti numbers.  If $\Gamma$ is a finitely generated residually finite group, it is known that
\[RG(\Gamma) \geq \text{cost}(\Gamma)-1 \geq \beta_{1}^{(2)}(\Gamma) - \frac{1}{|\Gamma|}\]
where we use the standard convention that $\frac{1}{|\Gamma|} = 0$ if $\Gamma$ is infinite.  Abert and Nikolov \cite{AbertNikolov} proved the first part of the inequality and the second part was proved by Gaboriau \cite{GabInequality}.  It is not known whether or not the inequalities can be strict. The reader is referred to \cite{Luck} for more information on L$^{2}$-Betti numbers.  We will discuss the notion of cost in Section~\ref{cost} as it will be central to our calculation of the rank gradient of amalgamated free products and HNN extensions. 

Namely, we prove and use the following lower bound for cost: 

\begin{proposition}
\label{cost prop}
Let $\Gamma$ be a finitely generated group and $L$ be a subgroup of $\Gamma$.  Let $\{H_{n}\}$ be a lattice of finite index normal subgroups of $\Gamma$ such that $\bigcap H_{n} = 1$.  Let $\widehat{\Gamma}_{(H_{n})}$ be the profinite completion of $\Gamma$ with respect to $\{H_{n}\}$ and define $\widehat{L}_{(L \cap H_{n})}$ similarly.  Then $\ds Cost(L, \widehat{\Gamma}_{(H_{n})}) = Cost(L, \widehat{L}_{(L \cap H_{n})}).$
\end{proposition}

Gaboriau \cite{GaboriauGroupCost} proved a lower bound for the cost of amalgamated free products and HNN extensions of groups over amenable subgroups. Proposition~\ref{cost prop} and Gaboriau's results are used in a fundamental way to compute the rank gradient of amalgamated free products and HNN extensions. The results are provided below.

\begin{theorem}
\label{amal intro}
Let $\Gamma = \Gamma_{1} \ast_{A} \Gamma_{2}$ be finitely generated and residually finite with $A$ amenable. Let $\{H_{n}\}$ be a lattice of normal subgroups of finite index in $\Gamma$ such that $\bigcap H_{n} = 1$.  Then
\[RG(\Gamma, \{H_{n}\}) = RG(\Gamma_{1}, \{\Gamma_{1} \cap H_{n}\}) + RG(\Gamma_{2}, \{\Gamma_{2} \cap H_{n}\}) + \frac{1}{|A|}.\]
In particular, $RG(\Gamma) \geq RG(\Gamma_{1}) + RG(\Gamma_{2}) + \frac{1}{|A|}$.
\end{theorem}

Let $K$ be a finitely generated group with isomorphic subgroups $A \simeq \varphi(A)$.  We denote the corresponding HNN extension of $K$ by $K \ast_{A}= \langle K, t \mid t^{-1}At = \varphi(A) \rangle$.  

\begin{theorem}
\label{hnn intro}
Let $\Gamma = K \ast_{A}$ be finitely generated and residually finite with $A$ amenable. Let $\{H_{n}\}$ be a lattice of finite index normal subgroups with $\bigcap H_{n} = 1$. Then 
\[RG(\Gamma, \{H_{n}\}) = RG(K, \{K \cap H_{n}\}) + \frac{1}{|A|}.\]  
In particular $RG(\Gamma) \geq RG(K) + \frac{1}{|A|}.$  
\end{theorem}

There are two key facts used to prove the above theorems.  The first is $RG(\Gamma) = cost(\Gamma, \widehat{\Gamma})-1$, which was proved by Abert and Nikolov \cite{AbertNikolov}. The second is the following theorem proved by Abert, Jaikin-Zapirain, and Nikolov in \cite{AJN}. 

\begin{theorem}[Abert, Jaikin-Zapirain, Nikolov]
\label{RGamenable}
Finitely generated infinite amenable groups have rank gradient zero with respect to any normal chain with trivial intersection.
\end{theorem}
Lackenby first proved this result for finitely presented groups in \cite{LacExpanders}.  The analogous statement about the $p$-gradient of infinite amenable groups also holds \cite{PappasArbpGrad}. 

Since there is no corresponding relationship between $p$-gradient and cost, the analogous results for the $p$-gradient of amalgamated free products and HNN extensions are proved differently.  In fact, the $p$-gradient is much easier to compute since $d_{p}(\Gamma)$ is easier to bound than $d(\Gamma)$.  To compute the $p$-gradient for amalgamated free products and HNN extensions we use the Kurosh subgroup theorems for amalgamated free products and HNN extensions.

Since $RG(A, \{A \cap H_{n}\}) = \frac{-1}{|A|}$ for amenable groups, then Theorem~\ref{amal intro} can be rewritten as 
\[RG(\Gamma_{1} \ast_{A} \Gamma_{2}, \{H_{n}\}) = RG(\Gamma_{1}, \{\Gamma_{1} \cap H_{n}\}) + RG(\Gamma_{2}, \{\Gamma_{2} \cap H_{n}\}) - RG(A, \{A \cap H_{n}\}).\]
  
The above equation does not hold in general as shown in Example~\ref{RG amalgam fails} by amalgamating over a subgroup with a large rank gradient.  However, if $\Gamma$ is the free product of an HNN extension with amenable associated subgroup and an infinite cyclic group then $\Gamma$ can be written as a non-trivial amalgamated free product.  Using Theorem~\ref{RG of HNN} one can show that the above equation for amalgamated products holds for $\Gamma$ when considered as an amalgamated free product. This illustrates that the condition of an amenable amalgamated subgroup is sufficient but not necessary.

The results given here are similar to the analogous results for cost proved by Gaboriau \cite{GaboriauGroupCost} and L$^{2}$-Betti numbers proved by L\"{u}ck \cite {LuckAppendix}. L\"{u}ck proved the corresponding equality of Theorem~\ref{amal intro} for the first L$^{2}$-Betti number of amalgamated free products and his result only requires that the first L$^{2}$-Betti number of the amalgamated subgroup is zero.   

\begin{acknowledgments}
The author would like to thank his advisor, Mikhail Ershov, for his help with the present material and also with his helpful comments on earlier drafts of this paper.  The author would also like to thank the anonymous referee for suggesting an improvement to Proposition~\ref{restricted cost}.
\end{acknowledgments}

\section{Rank Gradient and $p$-Gradient of Free Products}
We begin the section by giving the precise definition of $p$-gradient.  The notion of the $p$-gradient of a group for a prime number $p$ is also referred to in the literature as the \textit{mod-$p$ rank gradient} or \textit{mod-$p$ homology gradient}.  The reader should be careful as some authors define $p$-gradient differently \cite{LarispLarge}.  

\begin{definition}
Let $p$ be a prime.  The \textit{$p$-gradient} of $\Gamma$ relative to a lattice $\{H_{n}\}$ of $p$-power index normal subgroups is defined as 
\[RG_{p}(\Gamma, \{H_{n}\}) = \inf_{n} \frac{d_{p}(H_{n})-1}{[\Gamma : H_{n}]}\]
where $d_p(H) = d\left(H / \frattini{H} \right)$.  One can also define the \textit{absolute $p$-gradient}, $RG_{p}(\Gamma)$, where the infimum is taken over all normal subgroups of $p$-power index. 
\end{definition}

\begin{remark}
\label{Lattice Chain}
To prove results about the rank gradient (analogously $p$-gradient) with respect to a lattice $\{H_{n}\}$ of normal subgroups of finite index in $\Gamma$, it is enough to prove the result for a descending chain of subgroups from the lattice. Specifically, one can use the chain: 
$H_{1} \geq H_{1} \cap H_{2} \geq H_{1} \cap H_{2} \cap H_{3} \geq \dots$ \hfill $\qed$
\end{remark}

We can prove a result similar to Theorem~\ref{RG of free prod rel chain} (Abert, Jaikin-Zapirain, and Nikolov) by following a similar method of proof. Namely we prove the analogous result for the absolute rank gradient and $p$-gradient of arbitrary finitely generated groups.  

\begin{theorem}
\label{RG of Free Products}
Let $\Gamma_{1}$ and $\Gamma_{2}$ be finitely generated groups.  Let $\Gamma = \Gamma_{1} \ast\Gamma_{2}$. Then $RG(\Gamma) = RG(\Gamma_{1}) + RG(\Gamma_{2}) +1$.
\end{theorem}

\begin{proof}
This is a simple reduction of the proof given by Abert, Jaikin-Zapirain, and Nikolov in \cite{AJN}.  One just needs to show that for every pair of  finite index subgroups $H_{1} \leq \Gamma_{1}$ and $H_{2} \leq\Gamma_{2}$, there exists a finite index subgroup $H \leq \Gamma$ such that $H_{i} = H \cap \Gamma_{i}$ for $i=1,2$.  Such an $H$ can be constructed by using the natural map  $\varphi: \Gamma \to \Gamma_{1} \times \Gamma_{2}$ and pulling back $H_{1} \times H_{2}$.
\end{proof}

The analogous result to Theorem~\ref{RG of Free Products} for $p$-gradient is now stated.  

\begin{theorem}
Let $\Gamma_{1}$ and $\Gamma_{2}$ be finitely generated groups and $p$ a prime number.  Let $\Gamma = \Gamma_{1} \ast\Gamma_{2}$. Then $RG_{p}(\Gamma) = RG_{p}(\Gamma_{1}) + RG_{p}(\Gamma_{2}) +1$.
\end{theorem}

\begin{proof}
The proof is identical to the proof of Theorem~\ref{RG of Free Products} by replacing ``subgroups" with ``normal subgroups" and ``finite index" with ``$p$-power index."  To complete the proof as in Abert, Jaikin-Zapirain, and Nikolov used in \cite{AJN}, one needs the following facts:
$(1)$ Let $\Gamma$ be a finitely generated group and $H$ a $p$-power index normal subgroup.  Then $d_{p}(H)-1 \leq (d_{p}(\Gamma)-1)[\Gamma:H]$. $(2)$ Let $A \ast B$ be the free product of two finitely generated groups.  Then $d_{p}(A \ast B) = d_{p}(A) + d_{p}(B)$. \end{proof}

\section{Cost of Restricted Actions}
\label{cost}
To get a lower bound for the rank gradient of amalgamated free products and HNN extensions over amenable subgroups we use the notion of cost.  The notion of cost was first introduced by Levitt \cite{Levitt} and for more information the reader is referred to \cite{GaboriauGroupCost, AbertNikolov, Levitt}.  The following exposition of cost closely follows \cite{AbertNikolov}.  Throughout this section, all measures are assumed to be normalized with respect to the compact space on which they are defined.  

Let $\Gamma$ be a countable group that acts on a standard Borel probability space $(X, \mu)$ by measure preserving Borel automorphisms.  Define the equivalence relation $E$ on $X$ by
\[ xEy \; \text{if there exists} \; \gamma \in \Gamma \; \text{with} \; y = \gamma x. \]
The relation $E$ is a Borel equivalence relation and every equivalence class is countable. Since $E$ is a subset of $X \times X$, we can consider $E$ as a graph on $X$.

\begin{definition}
A \textit{Borel subgraph} of $E$ is a directed graph on $X$ such that the edge set is a Borel subset of $E$.
\end{definition}

\begin{definition}
A subgraph $S$ of $E$ is said to \textit{span} $E$, if for any $(x,y) \in E$ with $x \neq y$
there exists a path from $x$ to $y$ in $S$, where a \textit{path} from $x$ to $y$ in $S$ is defined as a sequence $x_{0}, x_{1}, \dots, x_{k} \in X$ such that: $x_{0} =x, \; x_{k} = y$; and $(x_{i}, x_{i+1}) \in S$ or $(x_{i+1}, x_{i}) \in S \; (0 \leq i \leq k - 1)$. 
\end{definition}

\begin{definition}
$S$ is called a \textit{graphing} of $E$ if it is a Borel subgraph of $E$ that spans $E$.  
\end{definition}

The edge-measure of a Borel subgraph $S$ of $E$ is defined as 
\[ e(S) = \int_{x \in X} deg_{S}(x)\; d\mu \]
where $deg_{S}(x)$ is the number of edges in $S$ with initial vertex $x$: 
\[deg_{S} (x) = |\{y \in X \mid (x, y) \in S\} |. \]
Note that $e(S)$ may be infinite. 

\begin{definition}
Let $\Gamma$ be a countable group acting on a standard Borel probability space $X$ by measure preserving Borel automorphisms.  Let $E$ denote the equivalence relation of this action.  The \textit{cost} of $E$ is defined as 
\[ Cost(E) = Cost(\Gamma, X) = \inf e(S) \]
where the infimum is taken over all graphings $S$ of $E$.
\end{definition}

Abert and Nikolov proved the following connection between rank gradient and cost.  Their actual result \cite[Theorem 1]{AbertNikolov} is more general than the special case given below, but the following is all that will be needed here.  

\begin{theorem}[Abert and Nikolov] 
\label{RG Cost}
Let $\Gamma$ be a finitely generated residually finite group and $\{H_{n}\}$ be a lattice of normal subgroups of finite index such that $\bigcap H_{n} = 1$.  Then 
\[ RG(\Gamma, \{H_{n}\}) = Cost(E) - 1 = Cost(\Gamma, \widehat{\Gamma}_{(H_{n})}) - 1, \] 
where $E$ is the equivalence relation coming from the action of $\Gamma$ on $\widehat{\Gamma}_{(H_{n})}$ (profinite completion of $\Gamma$ with respect to the lattice of subgroups $\{H_{n}\}$) by left multiplication and $\widehat{\Gamma}_{(H_{n})}$ comes with its normalized Haar measure.  
\end{theorem}

As the above theorem indicates, we will be interested in a group acting on its profinite completion by left multiplication.  Since a profinite group is a compact topological group, the following theorem about invariant measures on homogenous spaces holds.  This theorem is a special case of \cite[Corollary B.1.7]{DeLaHarpePropertyT}, which states the result for locally compact groups under an additional assumption that is satisfied by all compact groups.  The result is needed to prove how the cost of the relation changes when restricting to a subspace.  The reader is referred to \cite{DeLaHarpePropertyT} or \cite{HaarIntegral} for more information about invariant measure on homogenous spaces and Haar integrals. 

\begin{theorem}
\label{HaarDoubleInt}
Let $G$ be a compact group and $H$ a closed subgroup of $G$.  Let $H \backslash G$ be the space of right cosets of $H$.  A nontrivial regular right invariant measure on $H \backslash G$ exists and such a measure $\mu_{H \backslash G}$ on $H \backslash G$ is unique up to multiplication by a positive constant.  The measure $\mu_{H \backslash G}$ satisfies the following condition:
\[ \int_{G} f(x) \; d\mu_{G}(x) = \int_{H \backslash G} \left( \int_{H} f(hx) \;d\mu_{H}(h) \right) \;d\mu_{H \backslash G}(Hx)\]
where $\mu_{G}$ and $\mu_{H}$ are the unique normalized Haar measures on $G$ and $H$ respectively.  
\end{theorem}

\begin{note}
The formula in Theorem~\ref{HaarDoubleInt} makes sense only when the function $\varphi:G \to \mathbb{C}$ given by $\varphi(x) = \int_{H} f(hx) \;d\mu_{H}(h)$ is constant on right cosets of $H$.  The fact that $\varphi$ is constant on right cosets of $H$ follows from the fact that $\mu_{H}$ is right invariant.  
\end{note}

The following lemma concerning profinite completions is elementary.  The proof follows from residual finiteness and \cite[Corollary 1.1.8]{Ribes}.

\begin{lemma}
\label{profinite completion inclusion}
Let $\Gamma$ be finitely generated and let $\{H_{n}\}$ be a lattice of normal subgroups of finite index in $\Gamma$.  Let $L$ be a subgroup of $\Gamma$.  Then $\widehat{L}_{(L \cap H_{n})}$ is isomorphic to a closed subgroup of $\widehat{\Gamma}_{(H_{n})}$.  
\end{lemma}

The following lemma is used in order to determine the cost of a restricted action.   

\begin{lemma}
\label{graphings}
Let $\Gamma$ be a finitely generated residually finite group and $\{H_{n}\}$ be a lattice of normal subgroups of finite index in $\Gamma$ such that $\bigcap H_{n} = 1$. Let $L$ be a subgroup of $\Gamma$ acting on $\widehat{\Gamma}_{(H_{n})}$ by left multiplication and denote the equivalence relation by $E_{L}^{\widehat{\Gamma}_{(H_{n})}}$.  Let $S$ be a graphing of $E_{L}^{\widehat{\Gamma}_{(H_{n})}}$. Let $\{\bar{g}\}$ denote a set of right coset representatives for $\widehat{L}_{(L \cap H_{n})}$ in $\widehat{\Gamma}_{(H_{n})}$. 
For any $\bar{g}$, let 
\[S_{\bar{g}} = \{(x,y) \in \widehat{L}_{(L \cap H_{n})} \times \widehat{L}_{(L \cap H_{n})} \mid (x \bar{g}, y \bar{g}) \in S \}.\]
Then 
\begin{enumerate} 
\item $S_{\bar{g}}$ is a graphing for $E_{L}^{\widehat{L}_{(L \cap H_{n})}}$.
\item Every graphing $S'$ of $E_{L}^{\widehat{L}_{(L \cap H_{n})}}$ is obtained as $S' = S_{t}$ for some graphing $S$ of $E_{L}^{\widehat{\Gamma}_{(H_{n})}}$ and a specific right transversal $T = \{t\}$ for $\widehat{L}_{(L \cap H_{n})}$ in $\widehat{\Gamma}_{(H_{n})}$.
\end{enumerate}
\end{lemma}

\begin{proof}
For notational simplicity, let $\widehat{\Gamma} = \widehat{\Gamma}_{(H_{n})}$ and $\widehat{L} = \widehat{L}_{(L \cap H_{n})}$.
\begin{enumerate}
\item Note that $(x,y) \in S_{\bar{g}}$ if and only if $(x \bar{g}, y \bar{g}) \in S$.  Also, by Lemma~\ref{profinite completion inclusion} it follows that $\widehat{L}$ is a closed subgroup of $\widehat{\Gamma}$. \\
\underline{Spanning:} 
We need to show that $S_{\bar{g}}$ spans $E_{L}^{\widehat{L}}$.   Let $(x,y) \in E_{L}^{\widehat{L}}$. Then there exists $\alpha \in L$ such that $\alpha x = y$ which implies $\alpha x \bar{g} = y \bar{g}$ and therefore $(x\bar{g}, y\bar{g}) \in E_{L}^{\widehat{\Gamma}}$.  Since $S$ is a graphing of $E_{L}^{\widehat{\Gamma}}$, then $S$ spans $E_{L}^{\widehat{\Gamma}}$.  Therefore, there exists a path from $x\bar{g}$ to $y\bar{g}$ in $S$, call it 
\[z_{0}, z_{1}, \dots z_{k}.\]  
By definition $z_{0} = x\bar{g}$, $z_{k} = y\bar{g}$, and $(z_{i}, z_{i+1})$ or $(z_{i+1}, z_{i}) \in S$ for $0 \leq i \leq k-1$.  Let $z_{i}' = z_{i}\bar{g}^{-1}$.  Then the path in $S$ from $x\bar{g}$ to $y\bar{g}$ is now 
\[z_{0}'\bar{g}, z_{1}'\bar{g}, \dots, z_{k}'\bar{g}.\]  
It follows that $z_{0}'\bar{g} = x\bar{g}$, $z_{k}'\bar{g} = y\bar{g}$, and $(z_{i}'\bar{g}, z_{i+1}'\bar{g})$ or $(z_{i+1}'\bar{g}, z_{i}'\bar{g}) \in S$ for $0 \leq i \leq k-1$.  Thus, $(z_{i}', z_{i+1}')$ or $(z_{i+1}', z_{i}') \in S_{\bar{g}}$ for $0 \leq i \leq k-1$.  Therefore there is a path in $S_{\bar{g}}$ from $x$ to $y$ and thus $S_{\bar{g}}$ spans $E_{L}^{\widehat{L}}$. \\
\underline{Borel Subgraph:} 
We need to show that the edge set of $S_{\bar{g}}$ is a Borel subset of $E_{L}^{\widehat{L}}$.  Let $\pi_{\bar{g}} : \widehat{L} \times \widehat{L} \to \widehat{\Gamma} \times \widehat{\Gamma}$ be given by $\pi_{\bar{g}}(x,y) = (x\bar{g}, y\bar{g})$.  Note that $\pi_{\bar{g}}$ is injective since $\widehat{L} \leq \widehat{\Gamma}$.  Since these spaces are topological groups, multiplication is a continuous map and so $\pi_{\bar{g}}$ is continuous.  By definition, $S_{\bar{g}} = \pi_{\bar{g}}^{-1}(S)$.  By continuity of $\pi_{\bar{g}}$ and the fact that $S$ is a Borel subgraph of $E_{L}^{\widehat{\Gamma}}$, it follows that $S_{\bar{g}}$ is a Borel subgraph of $E_{L}^{\widehat{L}}$.

\item 
Let $\phi: \widehat{\Gamma} \to \widehat{\Gamma}/\widehat{L}$ be the natural projection.  By \cite[Proposition I.1.1]{SerreCoho}, it follows that  $\phi$ admits a continuous section $s: \widehat{\Gamma}/\widehat{L} \to \widehat{\Gamma}$. Then set $T=Im(s)$ is a (closed) right transversal for $\widehat{L}$ in $\widehat{\Gamma}$. Using this transversal, every element $g \in\widehat{\Gamma}$ can be uniquely written as $g = \ell t$ for some $t \in T$ and $\ell = g (s (\phi(g)))^{-1} \in \widehat{L}$.  The map $\psi: \widehat{\Gamma} \to \widehat{L}$ given by $g \to \ell$ is continuous by construction.
%
%

Let $S'$ be a graphing of $E_{L}^{\widehat{L}}$.  For every $t \in T$, we can form a graphing $S't$ of the right coset $\widehat{L}t$ in $\widehat{\Gamma}$. Consider $\ds S = \bigcup_{t \in T} S't$.  Then $S = \psi^{-1}(S')$ and since $\psi$ is continuous, it follows that $S$ is a graphing of $E_{L}^{\widehat{\Gamma}}$.  Furthermore, $S' = S_{t_{0}}$.  
\end{enumerate}
\end{proof}

Using the above theorem and lemma we can now prove the following result about the cost of a restricted action.  

\begin{proposition}
\label{restricted cost}
Let $\Gamma$ be a finitely generated group and $L$ be a subgroup.  Let $\{H_{n}\}$ be a lattice of finite index normal subgroups of $\Gamma$ such that $\bigcap H_{n} = 1$.  Let $\widehat{\Gamma}_{(H_{n})}$ be the profinite completion of $\Gamma$ with respect to $\{H_{n}\}$ and define $\widehat{L}_{(L \cap H_{n})}$ similarly.  Then $\ds Cost(L, \widehat{\Gamma}_{(H_{n})}) = Cost(L, \widehat{L}_{(L \cap H_{n})}).$
\end{proposition}

\begin{proof}
Let 
\[deg_{R}^{X} (x) = |\{y \in X \mid (x, y) \in R\} | \]
for any graphing $R$ on $E_{G}^{X}$, where $G$ is a group acting on the space $X$.  Let $\widehat{\Gamma} = \widehat{\Gamma}_{(H_{n})}$ and let $\widehat{L} = \widehat{L}_{(L \cap H_{n})}$.  By Lemma~\ref{profinite completion inclusion}, $\widehat{L}$ is a closed subgroup of $\widehat{\Gamma}$.

Let $T$ be the right transversal of $\widehat{L}$ in $\widehat{\Gamma}$ as constructed in Lemma~\ref{graphings}.  We know that if $S$ is a graphing of $E_{L}^{\widehat{\Gamma}}$, then $S_{t}$ is a graphing of $E_{L}^{\widehat{L}}$ for every $t \in T$.  For $g \in \widehat{\Gamma}$,  there is a map $g \to (\ell_{g}, \widehat{L}t_{g}) \in \widehat{L} \times \widehat{L} \backslash \widehat{\Gamma}$ where $\ell_{g}t_{g}= g$.  

For $(\ell, \widehat{L} t_{g}) \in \widehat{L} \times \widehat{L} \backslash \widehat{\Gamma}$ set
\[deg_{S}^{\widehat{\Gamma}}(\ell, \widehat{L} t_{g}) = deg_{S}^{\widehat{\Gamma}}(\ell t_{g}) = |\{x \in \widehat{\Gamma} \mid (\ell t_{g},x) \in S\}|. \]
Fix $\widehat{L} t_{g} \in \widehat{L} \backslash \widehat{\Gamma}$.  Then
\begin{align*}
deg_{S}^{\widehat{\Gamma}}(\ell, \widehat{L} t_{g}) 
&= |\{x \in \widehat{\Gamma} \mid (\ell t_{g},x) \in S\}| = |\{z \in \widehat{\Gamma} \mid (\ell t_{g}, z t_{g}) \in S\}| \\
(\ast) \hspace{15pt}
&= |\{y \in \widehat{L} \mid (\ell t_{g}, y t_{g}) \in S\}| = |\{y \in \widehat{L} \mid (\ell, y) \in S_{t_{g}}\}| \\
&= deg_{S_{t_{g}}}^{\widehat{L}}(\ell). 
\end{align*}
The equality ($\ast$) is given by the following:  Since $\widehat{L} \leq \widehat{\Gamma}$ it is clear that 
\[ \{y \in \widehat{L} \mid (\ell t_{g}, y t_{g}) \in S\} \subseteq \{z \in \widehat{\Gamma} \mid (\ell t_{g}, z t_{g}) \in S\} \]
and therefore we have the inequality $\geq$. Let $z \in \widehat{\Gamma}$ with $(\ell t_{g}, z t_{g}) \in S$.  Then  $\ell t_{g}, z t_{g} \in E_{L}^{\widehat{\Gamma}}$ and thus there is an $\alpha \in L$ such that $z t_{g} = \alpha \ell t_{g}$.  Thus $z = \alpha \ell \in \widehat{L}$ since $L \subset \widehat{L}$ by assumption.  The inequality $\leq$ follows.  Thus, for all $\widehat{L}t_{g} \in \widehat{L} \backslash \widehat{\Gamma}$ we have $deg_{S}^{\widehat{\Gamma}}(\ell, \widehat{L}t_{g}) = deg_{S_{t_{g}}}^{\widehat{L}}(\ell)$.  

Let $\mu_{\widehat{\Gamma}}$ and $\mu_{\widehat{L}}$ be the unique normalized Haar measures on $\widehat{\Gamma}$ and $\widehat{L}$ respectively.  By Lemma~\ref{profinite completion inclusion} it follows that $\widehat{L}$ is a closed subgroup of $\widehat{\Gamma}$ and therefore,
\begin{align*}
Cost(L, \widehat{\Gamma}) 
&= \inf_{\substack{S \; \text{graphing} \\ \text{of} \; E_{L}^{\widehat{\Gamma}}}} \int_{\widehat{\Gamma}} deg_{S}^{\widehat{\Gamma}}(g)\; d\mu_{\widehat{\Gamma}}(g) \\
\text{\scriptsize{by Theorem~\ref{HaarDoubleInt}}} \hspace{15pt}
&= \inf_{S} \int_{\widehat{L} \backslash \widehat{\Gamma}} \int_{\widehat{L}} deg_{S}^{\widehat{\Gamma}}(\ell, \widehat{L} t_{g})\; d\mu_{\widehat{L}}(\ell) \; d\mu_{\widehat{L} \backslash \widehat{\Gamma}}(\widehat{L} t_{g}) \\
\text{\scriptsize{by above}} \hspace{15pt}
&= \inf_{S} \int_{\widehat{L} \backslash \widehat{\Gamma}} \int_{\widehat{L}} deg_{S_{t_{g}}}^{\widehat{L}}(\ell) \; d\mu_{\widehat{L}}(\ell) \; d\mu_{\widehat{L} \backslash \widehat{\Gamma}}(\widehat{L} t_{g}) \\
(\ast \ast) \hspace{15pt} 
&= \inf_{S} \int_{\widehat{L} \backslash \widehat{\Gamma}} Cost(L, \widehat{L})\; d\mu_{\widehat{L} \backslash \widehat{\Gamma}}(\widehat{L} t_{g}) \\
&= Cost(L, \widehat{L}) \; \mu_{\widehat{L} \backslash \widehat{\Gamma}}(\widehat{L} \backslash \widehat{\Gamma}) \\
&= Cost(L, \widehat{L}).
\end{align*}
The equality $(\ast \ast)$ is as follows: Recall the definition of cost
\[Cost(L, \widehat{L}) 
= \inf_{\substack{S' \; \text{graphing} \\ \text{of} \; E_{L}^{\widehat{L}}}} \int_{\widehat{L}} deg_{S'}^{\widehat{L}}(\ell)\; d\mu_{1}(\ell). \]
The inequality $\geq$ is by the definition of cost and the inequality $\leq$ follows by definition and Proposition~\ref{graphings}.2.
\end{proof}

\section{Rank Gradient of Amalgamated Free Products and HNN Extensions}

\subsection{Rank Gradient of Amalgamated Free Products}
\label{RG Amals}
Let $\Gamma = \Gamma_{1} \ast_{A} \Gamma_{2}$ be residually finite and assume $A$ is amenable.  Let $\{H_{n}\}$ be a lattice of normal subgroups of finite index in $\Gamma$ such that $\bigcap H_{n} = 1$.  The action of $\Gamma$ on the boundary of the coset tree $\partial T(\Gamma, \{H_{n}\})$ is the action of $\Gamma$ by left multiplication on its profinite completion with respect to the lattice $\{H_{n}\}$ with normalized Haar measure.  For notational simplicity denote this completion and measure by $\widehat{\Gamma} = \widehat{\Gamma}_{(H_{n})}$ and $\mu$ respectively.  Since for $i=1, 2$,  $\{\Gamma_{i} \cap H_{n}\}$ is a lattice of finite index normal subgroups of $\Gamma_{i}$ with trivial intersection, then we have the completions $\widehat{\Gamma_{i}} = \widehat{\Gamma_{i}}_{(\Gamma_{i} \cap H_{n})}$ with measures $\mu_{i}$.  Similarly define $\widehat{A} = \widehat{A}_{(A \cap H_{n})}$ and $\mu_{A}$. Note that these completions are all profinite groups and thus are compact Hausdorff topological groups.  By Lemma~\ref{profinite completion inclusion} it follows that $\widehat{\Gamma_{i}} \leq \widehat{\Gamma}$.  

\begin{theorem}
\label{RG amalgams}
Let $\Gamma = \Gamma_{1} \ast_{A} \Gamma_{2}$ be finitely generated and residually finite with $A$ amenable.  Let $\{H_{n}\}$ be a lattice of normal subgroups of finite index in $\Gamma$ such that $\bigcap H_{n} = 1$.  Then
\begin{equation}
RG(\Gamma, \{H_{n}\}) = RG(\Gamma_{1}, \{\Gamma_{1} \cap H_{n}\}) + RG(\Gamma_{2}, \{\Gamma_{2} \cap H_{n}\}) + \frac{1}{|A|}.     \label{rg amalgam eq}
\end{equation}
In particular, $RG(\Gamma) \geq RG(\Gamma_{1}) + RG(\Gamma_{2}) + \frac{1}{|A|}.$
\end{theorem}

\begin{note}
This theorem was independently proved by Kar and Nikolov \cite[Proposition 2.2]{KarNikolov} in the case of amalgamation over a finite subgroup using Bass-Serre theory.  We will thus only show the case where $A$ is infinite amenable.  
\end{note}

\begin{proof}
Since $A$ is infinite we only need to show that 
\[RG(\Gamma, \{H_{n}\}) = RG(\Gamma_{1}, \{\Gamma_{1} \cap H_{n}\}) + RG(\Gamma_{2}, \{\Gamma_{2} \cap H_{n}\}).\]  

To simplify notation let $\widehat{\Gamma} = \widehat{\Gamma}_{(H_{n})}$ and $\widehat{\Gamma_{i}} = \widehat{\Gamma_{i}}_{(\Gamma_{i} \cap H_{n})}$ for $i=1,2$ and let $E = E_{\Gamma}^{\widehat{\Gamma}}$, $E|_{\Gamma_{i}} = E_{\Gamma_{i}}^{\widehat{\Gamma}}$, and $E|_{A} = E_{A}^{\widehat{\Gamma}}$.  Any action of an infinite amenable group on a Borel probability space is hyperfinite \cite{OrnsteinWeiss} and thus the cost is equal to 1; see \cite{GaboriauGroupCost}.  Therefore, $Cost(E|_{A})=1$.  

Theorem~\ref{RG Cost} states 
\[RG(\Gamma, \{H_{n}\}) = Cost(\Gamma, \widehat{\Gamma}) - 1 = Cost(E) - 1\]
and by \cite[Example 4.8]{GaboriauGroupCost}, it follows that $E = E|_{\Gamma_{1}} \ast_{E|_{A}} E|_{\Gamma_{2}}$.  
Since $E|_{A}$ is hyperfinite, then by \cite[Theorem 4.15]{GaboriauGroupCost}
\[Cost(E|_{\Gamma_{1}} \ast_{E|_{A}} E|_{\Gamma_{2}}) - 1 = Cost(E|_{\Gamma_{1}}) + Cost(E|_{\Gamma_{2}}) - Cost(E|_{A}) - 1.\]
Thus, 
\begin{align*}
RG(\Gamma, \{H_{n}\}) 
&= Cost(E|_{\Gamma_{1}} \ast_{E|_{A}} E|_{\Gamma_{2}}) - 1 \\
&= Cost(E|_{\Gamma_{1}}) + Cost(E|_{\Gamma_{2}}) - Cost(E|_{A}) - 1 \\
&= \left(Cost(E|_{\Gamma_{1}})-1\right) + \left(Cost(E|_{\Gamma_{2}}) -1\right) \\
&= \left(Cost(\Gamma_{1}, \widehat{\Gamma})-1\right) + \left(Cost(\Gamma_{2}, \widehat{\Gamma}) -1\right) \\
\text{\scriptsize by Prop~\ref{restricted cost}} \hspace*{15 pt}
&= \left(Cost(\Gamma_{1}, \widehat{\Gamma_{1}} )-1\right) + \left(Cost(\Gamma_{2}, \widehat{\Gamma_{2}} ) -1\right) \\
\text{\scriptsize by Theorem~\ref{RG Cost}} \hspace*{15 pt}
&=RG(\Gamma_{1}, \{\Gamma_{1} \cap H_{n}\}) + RG(\Gamma_{2}, \{\Gamma_{2} \cap H_{n}\}).
\end{align*}
Therefore, $RG(\Gamma, \{H_{n}\}) = RG(\Gamma_{1}, \{\Gamma_{1} \cap H_{n}\}) + RG(\Gamma_{2}, \{\Gamma_{2} \cap H_{n}\})$.  

The fact that $RG(\Gamma) \geq RG(\Gamma_{1}) + RG(\Gamma_{2}) + \frac{1}{|A|}$ follows by applying (\ref{rg amalgam eq}) to the lattice of all subgroups of finite index in $\Gamma$ and the definition of rank gradient.  
\end{proof}

\subsection{Rank Gradient of HNN Extensions}
\label{Rank Grad HNN}
Let $K$ be a finitely generated group with isomorphic subgroups $A \simeq \varphi(A)$.  We denote the associated HNN extension of $K$ by $K \ast_{A}= \langle K, t \mid t^{-1}At = \varphi(A) \rangle$.  Let $\{H_{n}\}$ be a lattice of finite index normal subgroups in $\Gamma = K \ast_{A}$ with $\bigcap H_{n} = 1$. Let $\widehat{\Gamma}_{(H_{n})}$ be the profinite completion of $\Gamma$ with respect to $\{H_{n}\}$ and let $\mu$ denote the unique normalized Haar measure on $\widehat{\Gamma}_{(H_{n})}$.  Define $\widehat{K}_{(K \cap H_{n})}$ and $\widehat{A}_{(A \cap H_{n})}$ similarly.

\begin{theorem}
\label{RG of HNN}
Let $\Gamma = K \ast_{A} = \langle K, t \mid t^{-1}At = B \rangle$ be finitely generated and residually finite with $A$ amenable.  Let $\{H_{n}\}$ be a lattice of finite index normal subgroups with $\bigcap H_{n} = 1$. Then 
\begin{equation}
RG(\Gamma, \{H_{n}\}) = RG(K, \{K \cap H_{n}\}) + \frac{1}{|A|}.    \label{rg hnn eq}  
\end{equation}
In particular, $RG(\Gamma) \geq RG(K) + \frac{1}{|A|}.$  
\end{theorem}

\begin{proof}
By Remark~\ref{Lattice Chain}, it is enough to prove the result assuming that $\{H_{n}\}$ is a descending chain. 

The result is proved analogously to Theorem~\ref{RG amalgams} using Gaboriau's results about the cost of HNN extensions \cite[Definition 4.20, Example 4.21, and Corollary 4.25]{GaboriauGroupCost}). 
\end{proof}

Recall since $RG(A, \{A \cap H_{n}\}) = \frac{-1}{|A|}$ for amenable groups, equation (\ref{rg amalgam eq}) can be written as 
\[RG(\Gamma_{1} \ast_{A} \Gamma_{2}, \{H_{n}\}) = RG(\Gamma_{1}, \{\Gamma_{1} \cap H_{n}\}) + RG(\Gamma_{2}, \{\Gamma_{2} \cap H_{n}\}) - RG(A, \{A \cap H_{n}\}).\]

The following example shows that the equation for amalgamated free products does not hold in general.

\begin{example}
\label{RG amalgam fails}
Let $\Gamma_{1} = F_{r} \times \mathbb{Z} / 2\mathbb{Z}, \Gamma_{2} = F_{r} \times \mathbb{Z} / 3\mathbb{Z}$, and let $A = F_{r}$.  Then $A$ is finite index in both $\Gamma_{1}$ and $\Gamma_{2}$ which implies
\begin{align*}
RG(\Gamma_{1}) +RG(\Gamma_{2}) - RG(A) 
&= \frac{RG(A)}{[\Gamma_{1}: A]} +\frac{RG(A)}{[\Gamma_{2}:A]} - RG(A) \\
&= \frac{r-1}{2} + \frac{r-1}{3} - (r-1) = -\frac{1}{6}(r-1).
\end{align*}
If we let $r=6k+1$, then $RG(\Gamma_{1}) +RG(\Gamma_{2}) - RG(A) = -k$ for any $k \in \mathbb{N}$.  However, for any finitely generated group $\Gamma$, we know  $RG(\Gamma) \geq -1$.  Therefore, $RG(\Gamma) \neq RG(\Gamma_{1}) +RG(\Gamma_{2}) - RG(A)$ in this case.  \qed
\end{example}

\section{$p$-Gradient of Amalgamated Free Products and HNN Extensions}
\label{p grad sec}
Computing the $p$-gradient for amalgamated free products and HNN extensions over amenable subgroups is easier than the rank gradient case since $d_{p}(G) = d(G/[G,G]G^{p})$ is easier to compute than $d(G)$ for any group $G$.  Specifically, one just needs to apply the Kurosh subgroup theorems for amalgamated free products and HNN groups \cite{Cohen}.  For our purposes we are only interested in applying the theorem to normal subgroups of finite index.  In this case we can state the theorems as follows:  

\underline{Amalgams:} Every normal subgroup $H$ of finite index in the amalgamated free product $\Gamma= \Gamma_{1} \ast_{A} \Gamma_{2}$ is an HNN group with base subgroup $L$ and $n = |H \backslash \Gamma / A| - |H \backslash \Gamma / \Gamma_{1}| - |H \backslash \Gamma / \Gamma_{2}|+ 1$ free generators with each associated subgroup being isomorphic to $A \cap H$.  Specifically, 
\[H = \langle L , t_{1}, \dots, t_{n} \mid t_{i}(A \cap H)t_{i}^{-1} = \varphi_{i}(A) \cap H \rangle \]
where the $\varphi_{i}$ are appropriate embeddings from $A$ to $L$.    

Further, $L$ is an amalgamated free product of $|H \backslash \Gamma / \Gamma_{1}|$  groups that are isomorphic to $\Gamma_{1} \cap H$ and $|H \backslash \Gamma / \Gamma_{2}|$ groups that are isomorphic to $\Gamma_{2} \cap H$ with at most $|H \backslash \Gamma / \Gamma_{1}| + |H \backslash \Gamma / \Gamma_{2}| - 1$ amalgamations each of which is isomorphic to $A \cap H$.  

\underline{HNN Extensions:} Every normal subgroup $H$ of finite index in the HNN extension $\Gamma= \langle K, t \mid tAt^{-1} = \varphi(A)\rangle$ is an HNN group with base subgroup $L$ and $n = |H \backslash \Gamma / A| - |H \backslash \Gamma / K|+ 1$ free generators with each associated subgroup being isomorphic to $A \cap H$.  Specifically, 
\[H = \langle L , t_{1}, \dots, t_{n} \mid t_{i}(A \cap H)t_{i}^{-1} = \varphi_{i}(A) \cap H \rangle \]
where the $\varphi_{i}$ are appropriate embeddings from $A$ to $L$.    

Further, $L$ is an amalgamated free product of $|H \backslash \Gamma / K|$  groups that are isomorphic to $K \cap H$ with at most $|H \backslash \Gamma / K| - 1$ amalgamations each of which is isomorphic to $A \cap H$.  

In either case, one can compute $d_{p}(H)$ by thinking of $H/[H,H]H^{p}$ as a vector space over $\mathbb{F}_{p}$ and thus $d_{p}(H)$ is the dimension of $H/[H,H]H^{p}$ over $\mathbb{F}_{p}$.  To compute $d_{p}(H)$ the following bounds for amalgamated free products and HNN extensions are used. 

\begin{proposition}
If $\Gamma_{1} \ast_{A} \Gamma_{2} = \langle \Gamma_{1}, \Gamma_{2} \mid A = \varphi(A) \rangle$ is an amalgamated free product, then 
\[d_{p}(\Gamma_{1}) + d_{p}(\Gamma_{2}) - d_{p}(A) \leq d_{p}(\Gamma_{1} \ast_{A} \Gamma_{2}) \leq d_{p}(\Gamma_{1}) + d_{p}(\Gamma_{2}).\]
If $K \ast_{A} =\langle K, t \mid tAt^{-1} = \varphi(A) \rangle$ is an HNN extension, then 
\[d_{p}(K) - d_{p}(A) + 1 \leq d_{p}(K \ast_{A}) \leq d_{p}(K) + 1.\]
\end{proposition}

Combining this proposition and the presentation of $H$ given by the Kurosh subgroup theorem allows one to compute $RG_{p}(\Gamma)$ directly by bounding $d_{p}(H)$ for any normal subgroup of $p$-power index in $\Gamma$.  The equations for the $p$-gradient of a free product amalgamated over an amenable subgroup and an HNN extension with amenable associated subgroup are analogous to the results of Theorem~\ref{RG amalgams} and Theorem~\ref{RG of HNN} respectively.  One just needs to replace ``$RG$" with ``$RG_{p}$", `residually finite" with ``residually-$p$", and ``finite index" with ``$p$-power index." 

\begin{theorem}
\label{pgradient amalgams}
Let $\Gamma = \Gamma_{1} \ast_{A} \Gamma_{2}$ be finitely generated and residually-$p$ with $A$ amenable. Let $\{H_{n}\}$ be a lattice of normal subgroups of $p$-power index in $\Gamma$ such that $\bigcap H_{n} = 1$.  Then
\begin{equation}
RG_{p}(\Gamma, \{H_{n}\}) = RG_{p}(\Gamma_{1}, \{\Gamma_{1} \cap H_{n}\}) + RG_{p}(\Gamma_{2}, \{\Gamma_{2} \cap H_{n}\})  + \frac{1}{|A|}. 
\end{equation}
In particular, $RG_{p}(\Gamma) \geq RG_{p}(\Gamma_{1}) + RG_{p}(\Gamma_{2}) + \frac{1}{|A|}$. 
\end{theorem}

\begin{theorem}
Let $\Gamma = K \ast_{A} = \langle K, t \mid tAt^{-1} = B \rangle$ be finitely generated and residually-$p$ with $A$ amenable.  Let $\{H_{n}\}$ be a lattice of normal subgroups of $p$-power index in $\Gamma$ such that $\bigcap H_{n} = 1$.  Then
\begin{equation}
RG_{p}(\Gamma,\{H_{n}\}) = RG_{p}(K, \{K \cap H_{n}\}) + \frac{1}{|A|}.
\end{equation}
In particular, $RG_{p}(\Gamma) \geq RG_{p}(K) + \frac{1}{|A|}$.  
\end{theorem}

\bibliographystyle{amsplain}
\bibliography{amalgamscite}

\end{document}